\documentclass[pra,onecolumn,superscriptaddress,10pt,showpacs,floatfix]{revtex4-1}
\usepackage[pdftex,plainpages=false,colorlinks=true,citecolor=blue,linkcolor=blue,urlcolor=blue,filecolor=green,bookmarksopen=true]{hyperref}
\usepackage[utf8]{inputenc}
\usepackage{amsmath}
\usepackage[caption = false]{subfig}
\usepackage{graphicx,epstopdf}
\usepackage{blindtext}
\usepackage{lipsum}
\usepackage{amsfonts}
\usepackage{bbm}
\usepackage{amssymb}
\usepackage{enumerate}
\usepackage{color}
\usepackage{latexsym}
\usepackage{times,txfonts}

\begin{document}

\title{Determining the eigenvalues of a square matrix through known
	information of its submatrix}

\author{M. A. de Ponte}
\email{mickel.ponte@unesp.br}
\affiliation{Universidade Estadual Paulista (UNESP), Campus Experimental de Itapeva, 18409-010, Itapeva, S\~{a}o Paulo, Brazil}

\author{L. C. de Campos}
\affiliation{Universidade Estadual Paulista (UNESP), Campus Experimental de Itapeva, 18409-010, Itapeva, S\~{a}o Paulo, Brazil}

\begin{abstract}
In this paper we bring to light an unprecedented property of the eigenvalues
of a matrix $\mathbf{A}$ with the eigenvalues and eigenvectors of a
submatrix of $\mathbf{A}$. This property can be used, through the technique
developed here, to determine some of eigenvalues of $\mathbf{A}$ and, thus,
reduce the degree of the characteristic polynomial associated with this
matrix. This reduction can occur in two ways: when restrictions between some
elements of $\mathbf{A}$ are checked (or imposed) and/or when some
eigenvalues of the submatrix of $\mathbf{A}$ are degenerate. As an
application, we show how to obtain matrices of arbitrary size, not trivial,
that allow the algebraic calculation of eigenvalues and eigenvectors and how
to share a set of eigenvectors of a submatrix of $\mathbf{A}$ with matrix $%
\mathbf{A}$ itself, preserving norm, direction and sense.

\end{abstract}

\maketitle

\section{Introduction}

A lot of problems in physics, mathematics, and other sciences can be solved
through techniques developed in matrix algebra. In view of the wide
application of these concepts, whether from an algebraic or numerical point
of view, it is always important and welcome the development of new
alternative methods that facilitate or reduce the calculation time of the
eigenvalues associated with the matrices.

With respect to algebraic methods, it is known that for a $N\times N$ square
matrix its characteristic polynomial will be of the order $N$. Thus, by
Abel's theorem \cite{Abel} we are unable to find algebraically the
eigenvalues for $N>4$, except for particular cases in which the original
characteristic polynomial can be written as products of polynomials of order 
$\leq 4$ or redefined through a new variable that transforms the original
polynomial into another polynomial of order $\leq 4$. In the literature we
find some examples of these particular cases where $N$ may be arbitrary: $i)$
linear chain , $ii)$ circular chain , $iii)$ symmetric chain , and $iv)$
star or central chain \cite{1,2,3,4,5,6}. On the other hand, when a
numerical solution of an eigenvalue problem is desired, there are several
methods. Some of these methods, which we will briefly discuss below, can be
found in ref. \cite{Book}. In general, numerical methods are composed of two
steps: in the first step, the original matrix is transformed into a simpler
equivalent matrix, which has the same eigenvalues, and in the second step,
we use iterative methods for the determination of the eigenvalues. In order
to perform this first step there are the Given and Householder methods \cite{Book},
which transform the original matrix into a simpler tridiagonal matrix, while
the Jacob method transforms the original matrix into a diagonal matrix after
infinite iterative processes. The second step is generally performed by
expanding the determinant to obtain a sequence of polynomials that can be
calculated iteratively to obtain the eigenvalues. These methods are
generally used when you want to determine part of the eigenvalues. When we
want to determine all eigenvalues{}we must use Rutishauser's LR method and
Francis's QR method \cite{Book}.

In this work, our main objective is to present the scientific community, in
a very simple way, the technique that allows us to reduce one order of the
characteristic polynomial for each restriction satisfied or imposed. As we
will see later, these restriction are formed by relations between some
elements of the original matrix $\mathbf{A}$ with the eigenvectors
associated with a square matrix $\mathbf{B}$\textbf{,} obtained from $%
\mathbf{A}$ excluding rows and columns. Note, therefore, that this technique
is useful when it is possible to identify a matrix $\mathbf{B}$ of which its
eigenvalues and eigenvectors are already known.

Below we present some examples of how it is possible to reduce the degree of
the characteristic polynomial, already known in the literature, and then
present the idea of our method through a very simple example.

\section{How can the order of the characteristic polynomial be reduced?}

For any square matrix $\mathbf{A}$, its eigenvalues can be determined
through the roots of a characteristic polynomial defined by $p\left( \lambda
\right) =\det \left( \mathbf{A}-\lambda \mathbf{I}\right) =0$, where $%
\mathbf{I}$ is the identity matrix, which has the same order as the matrix $%
\mathbf{A}$. $p\left( \lambda \right) $ is a polynomial of degree $N$, which
is the size of matrix $\mathbf{A}$, and on certain conditions can be reduced
to a smaller order, depending on of the coefficients that multiply each
power of $\lambda $. Some examples of this reduction are shown below:

\begin{enumerate}
	\item By redefining a new variable as a function of $\lambda $, for example%
	\begin{equation*}
	\lambda ^{4}+a_{1}\lambda ^{2}+a_{2}=0\rightarrow \left\{ 
	\begin{array}{c}
	y\equiv \lambda ^{2} \\ 
	y^{2}+a_{1}y+a_{2}=0\qquad %
	\end{array}%
	\right.
	\end{equation*}
	is reduced to a 2nd order polynomial.
	
	\item Through knowledge of one of the roots of $p\left( \lambda \right) $,
	for instance 
	\begin{equation*}
	\lambda ^{4}+b_{1}\lambda ^{3}+b_{2}\lambda ^{2}+b_{3}\lambda =0\rightarrow
	\left\{ 
	\begin{array}{c}
	\lambda _{1}=0 \\ 
	\lambda ^{3}+b_{1}\lambda ^{2}+b_{2}\lambda +b_{3}=0\qquad%
	\end{array}%
	\right.
	\end{equation*}
	is reduced to a 3rd order polynomial.
\end{enumerate}

In order to make clear our method for reducing the degree of the
characteristic polynomial $p\left( \lambda \right) $, which we will present
in a more general way later, we would like to illustrate it through a simple
example of a $3\times 3$ matrix. Next, we perform the demonstration for the
general case of an arbitrary $N\times N$ square matrix.

\subsection{An illustrative example}

Consider an arbitrary $3\times 3$ square matrix $\mathbf{A}$ with complex
elements $a_{i,j}$ which does not need to be Hermitian, but which may be
diagonalizable. Since the eigenvalues do not change when we interchange two
lines $i\leftrightarrow j$ simultaneously with the permutation of two
columns $i\leftrightarrow j$, we can rearrange the elements of matrix $%
\mathbf{A}$ in the most convenient way and thus consider a submatrix $%
\mathbf{B}$, obtained from $\mathbf{A}$ through elimination, for example,
from the 3rd row and 3rd column, that is,%
\begin{equation*}
\mathbf{B}=\left( 
\begin{array}{cc}
a_{1,1} & a_{1,2} \\ 
a_{2,1} & a_{2,2}%
\end{array}%
\right) =\left( 
\begin{array}{cc}
b_{1,1} & b_{1,2} \\ 
b_{2,1} & b_{2,2}%
\end{array}%
\right) \text{.}
\end{equation*}%
The matrix $\mathbf{A}$ will then be defined by%
\begin{equation*}
\mathbf{A}=\left( 
\begin{array}{cc}
\mathbf{B} & 
\begin{array}{c}
a_{1,3} \\ 
a_{2,3}%
\end{array}
\\ 
\begin{array}{cc}
a_{3,1} & a_{3,2}%
\end{array}
& a_{3,3}%
\end{array}%
\right) =\left( 
\begin{array}{ccc}
b_{1,1} & b_{1,2} & a_{1,3} \\ 
b_{2,1} & b_{2,2} & a_{2,3} \\ 
a_{3,1} & a_{3,2} & a_{3,3}%
\end{array}%
\right) \text{.}
\end{equation*}

This choice of $\mathbf{B}$, through the rearrangement of $\mathbf{A}$
elements, is desirable once we know its eigenvalues and eigenvectors, as
will become clear later. Let us then call the eigenvalues of $\mathbf{B}$ by 
$\lambda _{1}$ and $\lambda _{2}$, while the respective eigenvectors are $%
\vec{\upsilon}_{1}=\left( \upsilon _{1;1},\upsilon _{1;2}\right) ^{\intercal
}$ and $\vec{\upsilon}_{2}=\left( \upsilon _{2;1},\upsilon _{2;2}\right)
^{\intercal }$.

Once the matrix $\mathbf{B}$ has been identified, we can now ``investigate"
how the knowledge of its eigenvalues and eigenvectors allows us to reduce
the degree of the characteristic polynomial $p\left( \lambda \right) $. An
initial attempt to perform this investigation is to analyze the necessary
conditions for an eigenvector $\vec{\vartheta}$, of matrix $\mathbf{A}$, to
be of the form%
\begin{equation}
\vec{\vartheta}\equiv \left( \vec{\upsilon}_{j},\mu \right) ^{\intercal
}=\left( \upsilon _{j;1},\upsilon _{j;2},\mu \right) ^{\intercal },
\label{exc}
\end{equation}%
where $\vec{\upsilon}_{j}=\left( \upsilon _{j;1},\upsilon _{j;2}\right)
^{\intercal }$ is one of the known eigenvectors of matrix $\mathbf{B}$. In
this section we assume that $i,j$ and $k$ are indices that can assume values 
$1$ and $2$.

If $\vec{\vartheta}$ ($\vec{\upsilon}$) is an eigenvector of matrix $\mathbf{%
	A}$ ($\mathbf{B}$), then we know that $\mathbf{A}\cdot \vec{\vartheta}%
=\Lambda \vec{\vartheta}$ ($\mathbf{B}\cdot \vec{\upsilon}=\lambda \vec{%
	\upsilon}$), where $\Lambda $ ($\lambda $) is an eigenvalue of matrix $%
\mathbf{A}$ ($\mathbf{B}$). We can therefore conclude that:%
\begin{equation}
\left( 
\begin{array}{ccc}
b_{1,1} & b_{1,2} & a_{1,3} \\ 
b_{2,1} & b_{2,2} & a_{2,3} \\ 
a_{3,1} & a_{3,2} & a_{3,3}%
\end{array}%
\right) \cdot \left( 
\begin{array}{c}
\upsilon _{j;1} \\ 
\upsilon _{j;2} \\ 
\mu%
\end{array}%
\right) =\Lambda \left( 
\begin{array}{c}
\upsilon _{j;1} \\ 
\upsilon _{j;2} \\ 
\mu%
\end{array}%
\right) \text{.}  \label{s1}
\end{equation}

Using the fact that $\mathbf{B}\cdot \vec{\upsilon}_{j}=\lambda _{j}\vec{%
	\upsilon}_{j}$, we can obtain from (\ref{s1}) the following system of
equations: 
\begin{subequations}
	\label{sub1}
	\begin{eqnarray}
	a_{13}\mu &=&\left( \Lambda -\lambda _{j}\right) \upsilon _{j;1}
	\label{sub1a} \\
	a_{23}\mu &=&\left( \Lambda -\lambda _{j}\right) \upsilon _{j;2}
	\label{sub1b} \\
	\left( \Lambda -a_{33}\right) \mu &=&a_{31}\upsilon _{j;1}+a_{32}\upsilon
	_{j;2}\text{.}  \label{sub1c}
	\end{eqnarray}
	
	From the equation (\ref{sub1a}) and (\ref{sub1b}) we realize that to exist a
	solution that satisfies (\ref{exc}) both equations must be linearly
	dependent (LD), that is, we must impose that 
\end{subequations}
\begin{equation}
a_{1,3}\upsilon _{j;2}=a_{2,3}\upsilon _{j;1}\text{,}  \label{V}
\end{equation}%
so that a connection is created between the eigenvector $\vec{\upsilon}_{j}$%
, of matrix $\mathbf{B}$, with the elements $a_{1,3}$ and $a_{2,3}$. From
the equation (\ref{V}) we see that the solution will be of the form%
\begin{equation*}
a_{1,3}=\alpha \upsilon _{j;1};\qquad a_{2,3}=\alpha \upsilon _{j;2},
\end{equation*}%
where $\alpha \in 
%TCIMACRO{\U{2102} }%
%BeginExpansion
\mathbb{C}
%EndExpansion
$ is a constant of proportionality. Note that this implies that the first
two elements of the 3rd column of matrix $\mathbf{A}$ are related to the
coefficients of the eigenvector $\vec{\upsilon}_{j}$, i.e.:%
\begin{equation}
\left( 
\begin{array}{c}
a_{1,3} \\ 
a_{2,3}%
\end{array}%
\right) =\alpha \left( 
\begin{array}{c}
\upsilon _{j;1} \\ 
\upsilon _{j;2}%
\end{array}%
\right) =\alpha \vec{\upsilon}_{j}\text{.}  \label{VV}
\end{equation}%
In addition to the equation (\ref{VV}), note that if we multiply (\ref{sub1c}%
) by $a_{2,3}$ on both sides of this equality and make use of equation (\ref%
{sub1b}), we obtain a second-degree equation for the eigenvalue $\Lambda $,
that is:%
\begin{equation}
\left( \Lambda -\lambda _{j}\right) \left( \Lambda -a_{3,3}\right) \upsilon
_{j;2}=a_{2,3}\left( a_{3,1}\upsilon _{j;1}+a_{3,2}\upsilon _{j;2}\right) 
\text{,}  \label{f1}
\end{equation}%
which can also be written with the help of (\ref{V}), in the form%
\begin{equation}
\left( \Lambda -\lambda _{j}\right) \left( \Lambda -a_{3,3}\right) \upsilon
_{j;1}=a_{1,3}\left( a_{3,1}\upsilon _{j;1}+a_{3,2}\upsilon _{j;2}\right) 
\text{.}  \label{f2}
\end{equation}%
If $\upsilon _{j;1}$ and/or $\upsilon _{j;2}$ are/is different from zero
(since both can not be null), using (\ref{V}) in (\ref{f2}) or (\ref{f1}),
we obtain%
\begin{equation}
\left( \Lambda -\lambda _{j}\right) \left( \Lambda -a_{3,3}\right)
=a_{1,3}a_{3,1}+a_{2,3}a_{3,2}\text{.}  \label{f3}
\end{equation}%
If the matrix $\mathbf{A}$ is Hermitian, that is, $a_{j,i}=a_{i,j}^{\ast }$,
the right side of equation (\ref{f3}) can be simplified using (\ref{VV}) for
the form 
\begin{equation}
\left( \Lambda -\lambda _{j}\right) \left( \Lambda -a_{3,3}\right)
=\left\vert \alpha \right\vert ^{2}\text{,}  \label{Auto}
\end{equation}%
since we are assuming that the set of eigenvectors $\left\{ \vec{\upsilon}%
_{j}\right\} $ are orthonormal, that is: $\vec{\upsilon}_{i}^{\dag }\cdot 
\vec{\upsilon}_{j}=\delta _{ij}$.

Note, therefore, that in the particular situation where the equality (\ref%
{VV}) is valid that for a given eigenvalue $\lambda _{j}$ and its respective
eigenvector $\vec{\upsilon}_{j}$ from matrix $\mathbf{B}$, we can determine
the eigenvalue $\Lambda $ of $\mathbf{A}$ through the equation (\ref{f3}).
Finally, we can obtain $\mu $ by the equation (\ref{sub1a}) or (\ref{sub1b}%
). In this way we see that this method provides $2$ eigenvalues for the
third-order square matrix, where the remaining root can be determined by
factoring the characteristic polynomial $p\left( \lambda \right) $
associated with the original matrix $\mathbf{A}$. However, the result for
some numerical examples has shown us that this remaining root is exactly the
other eigenvalue of $\mathbf{B}$, which we will call $\lambda _{k\neq j}$.
To confirm this result analytically, we verified that the characteristic
polynomial associated to matrix $\mathbf{A}$, using the Laplace development,
can be written in the form:%
\begin{eqnarray}
\det \left( \mathbf{A}-\Lambda \mathbf{I}\right)  &=&\left( a_{3,3}-\Lambda
\right) \det \left( \mathbf{B}-\Lambda \mathbf{I}\right)+a_{3,2}\left[
b_{2,1}a_{1,3}-a_{2,3}\left( b_{1,1}-\Lambda \right) \right]+a_{3,1}\left[ b_{1,2}a_{2,3}-a_{1,3}\left( b_{2,2}-\Lambda \right) \right]
\text{.}  \label{4}
\end{eqnarray}%
From the eigenvector equation for matrix $\mathbf{B}$ it can be shown that 
\begin{subequations}
	\label{Eq5}
	\begin{eqnarray}
	b_{1,2}\upsilon _{j;2}=\left( \lambda
	_{j}-b_{1,1}\right) \upsilon _{j;1}  \label{Eq5a} \\
	b_{2,1}\upsilon _{j;1}=\left( \lambda
	_{j}-b_{2,2}\right) \upsilon _{j;2}\text{,}  \label{Eq5b}
	\end{eqnarray}%
	so that by multiplying both sides of equalities in (\ref{Eq5}) by $\alpha $
	and using the relations in (\ref{VV}), we can rewrite the equation (\ref{4})
	in the form 
\end{subequations}
\begin{eqnarray}
\det \left( \mathbf{A}-\Lambda \mathbf{I}\right)  &=&\left( \Lambda +\lambda _{j}-b_{1,1}-b_{2,2}\right) \left(
a_{1,3}a_{3,1}+a_{2,3}a_{3,2}\right)+\left( a_{3,3}-\Lambda
\right) \det \left( \mathbf{B}-\Lambda \mathbf{I}\right) \text{.}  \label{e5}
\end{eqnarray}

If we recall that $Tr\mathbf{B}=b_{1,1}+b_{2,2}=\lambda
_{j}+\lambda _{k\neq j}$, then (\ref{e5}) can be rewritten as 
\begin{eqnarray}
\det \left( \mathbf{A}-\Lambda \mathbf{I}\right) &=&\left( \Lambda
-\lambda _{k\neq j}\right) \left( a_{1,3}a_{3,1}+a_{2,3}a_{3,2}\right)+\left( a_{3,3}-\Lambda
\right) \det \left( \mathbf{B}-\Lambda \mathbf{I}\right) \text{.}
\label{Cara}
\end{eqnarray}

This equation (\ref{Cara}) shows that for $\Lambda =\lambda _{k\neq j}$ we
obtain $\det \left( \mathbf{A}-\Lambda \mathbf{I}\right) =0$, since the term 
$\det \left( \mathbf{B}-\Lambda \mathbf{I}\right) $ also cancels, seen that $%
\Lambda =\lambda _{k\neq j}$ is also the eigenvalue of matrix $\mathbf{B}$.
This shows how to obtain the three eigenvalues of matrix $\mathbf{A}$. The
eigenvector associated with this eigenvalue, which will not be of the form (%
\ref{exc}) like the other eigenvectors, can be found through the solution of
the linear equations system:%
\begin{eqnarray}
b_{1,1}x+b_{1,2}y+a_{1,3}z &=&\lambda _{k\neq j}x  \notag \\
b_{2,1}x+b_{2,2}y+a_{2,3}z &=&\lambda _{k\neq j}y  \label{EL} \\
a_{3,1}x+a_{3,2}y+a_{3,3}z &=&\lambda _{k\neq j}z.  \notag
\end{eqnarray}%
A particular solution of this system (\ref{EL}) arises when the matrix $%
\mathbf{A}$ is ``almost Hermitian", that is, when $a_{j,i}=\beta
a_{i,j}^{\ast }$, irrespective of the elements $b_{i,j}$. Therefore, if we
use (\ref{VV}) we have that $a_{3,1}=\beta a_{1,3}^{\ast }=\beta \alpha
^{\ast }\upsilon _{j;1}^{\ast }$ and $a_{3,2}=\beta a_{2,3}^{\ast }=\beta
\alpha ^{\ast }\upsilon _{j;2}^{\ast }$, such that the amount 
\begin{equation}
a_{3,1}\upsilon _{k;1}+a_{3,2}\upsilon _{k;2}=\beta \alpha ^{\ast }\vec{%
	\upsilon}_{j}^{\dag }\cdot \vec{\upsilon}_{k}=0  \label{A}
\end{equation}%
when we consider $j\neq k$. Using (\ref{A}) we notice that the system (\ref%
{EL}) will have the following solution: $x=\upsilon _{k;1}$,$y=\upsilon
_{k;2}$, and $z=0$, where we then conclude that the eigenvector of matrix $%
\mathbf{B}$ is ``transferred" to matrix $\mathbf{A}$ in an ``intact" way, that
is, the vector continues to have the same modulus, direction and sense, and
only we need to add a new coordinate $z=0$ due to its extension from the
space $%
%TCIMACRO{\U{2102} }%
%BeginExpansion
\mathbb{C}
%EndExpansion
^{2}$ to the space $%
%TCIMACRO{\U{2102} }%
%BeginExpansion
\mathbb{C}
%EndExpansion
^{3}$.

The conclusion that $\lambda _{k\neq j}$ is an eigenvalue shared between
matrices $\mathbf{A}$ and $\mathbf{B}$ and that its eigenvector is the same
vector, in the ``almost Hermitian" situation, was a surprise for us. This led
us to consider the more general situation in which the third column
(excluding the element of principal diagonal) can be written as a linear
combination of all eigenvectors of matrix $\mathbf{B}$, as we shown below.

\section{One step further: Developing the technique in a completely way}

We saw above, through (\ref{VV}), that when the first two elements of the
3rd column of matrix $\mathbf{A}$ are proportional to one of the
eigenvectors of $\mathbf{B}$, say $\vec{\upsilon}_{j}$, that $2$ of the
eigenvalues of $\mathbf{A}$ can be found through a second-degree polynomial
as a function of the eigenvalue $\lambda _{j}$ of $\mathbf{B}$, whose
eigenvector is $\vec{\upsilon}_{j}$ (\ref{f3}). The third eigenvalue of $%
\mathbf{A}$ is exactly the eigenvalue of $\mathbf{B}$ associated with the
other eigenvector $\vec{\upsilon}_{k\neq j}$ of $\mathbf{B}$. We can now
rephrase these statements in another way: imagining that the first two
elements of the third column of matrix $\mathbf{A}$ (which we are admitting
that at least one of the elements is nonzero) compose a vector $\vec{\nu}$
in the space $%
%TCIMACRO{\U{2102} }%
%BeginExpansion
\mathbb{C}
%EndExpansion
^{2}$, we can always write it as a linear combination of the eigenvectors of 
$\mathbf{B}$, that is: $\vec{\nu}=\alpha \vec{\upsilon}_{1}+\beta \vec{%
	\upsilon}_{2}$, where $\vec{\upsilon}_{1}$ and $\vec{\upsilon}_{2}$ are
eigenvectors of $\mathbf{B}$ while $\alpha $ and $\beta $ are arbitrary
complex constants. In this sense, we can state that when $\alpha $ ($\beta $%
) is null, that the eigenvalue associated with $\vec{\upsilon}_{1}$ ($\vec{%
	\upsilon}_{2}$) will also be eigenvalue of $\mathbf{A}$ while the other
eigenvalues{}can be obtained through a second degree equation, instead of
3rd degree, which is written as a function of the eigenvalue associated with
the eigenvector $\vec{\upsilon}_{2}$ ($\vec{\upsilon}_{1}$), of which $\vec{%
	\nu}$ is a linear combination. At this point, we should be thinking: ``what
happens when both $\alpha $ and $\beta $ are non-zero?". The answer will be
given soon after, when we expose the technique in general.

\subsection{Reducing the degree of the characteristic polynomial}

Since the notation used in the general demonstration becomes dense and could
be difficult to understand, we have chosen to divide it into $3$ steps for
convenience.

\subsubsection{Step 1: Case in which $M=N+1$ and the elements $\left\{
	a_{i,M}\right\} $ are proportional to a single eigenvector of $\mathbf{B}$}

Consider a $M\times M$ square matrix $\mathbf{A}$ diagonalizable with
complex elements $a_{i,j}$. Let us assume that the $N$ first rows and
columns (initially fixing $N=M-1$) compose a square matrix $\mathbf{B}$,
defined by%
\begin{equation*}
\mathbf{B}=\left( 
\begin{array}{cccc}
b_{1,1} & b_{1,2} & \cdots  & b_{1,N} \\ 
b_{2,1} & b_{2,2} & \cdots  & b_{2,N} \\ 
\vdots  & \vdots  & \ddots  & \vdots  \\ 
b_{N,1} & b_{N,2} & \cdots  & b_{N,N}%
\end{array}%
\right) \text{,}
\end{equation*}%
such that the matrix $\mathbf{A}$ can be written in the form%
\begin{equation*}
\mathbf{A}=\left( 
\begin{array}{cc}
\mathbf{B} & 
\begin{array}{c}
a_{1,M} \\ 
a_{2,M} \\ 
\vdots  \\ 
a_{N,M}%
\end{array}
\\ 
\begin{array}{cccc}
a_{M,1} & a_{M,2} & \cdots  & a_{M,N}%
\end{array}
& a_{M,M}%
\end{array}%
\right) \text{.}
\end{equation*}%
Let us assume that we know the eigenvalues and eigenvectors of matrix $%
\mathbf{B}$, so that the $k$-th eigenvalue of $\mathbf{B}$ will be
represented by $\lambda _{k}$ and the respective eigenvector by $\vec{%
	\upsilon}_{k}=\left( \upsilon _{k;1},\upsilon _{k;2},\ldots ,\upsilon
_{k;N}\right) ^{\intercal }$. Assuming that one of the eigenvectors of $%
\mathbf{A}$ is of the form $\vec{\vartheta}_{k}=\left( \vec{\upsilon}%
_{k},\mu \right) ^{\intercal }$, we want to know what conditions this
hypothesis implies! (From this section onwards we will consider that the
indices $i,j$ and $k$ vary from $1$ to $N$.)

Before performing this analysis, note that if $\vec{\upsilon}_{k}$ is an
eigenvector of $\mathbf{B}$, which%
\begin{equation}
\sum\nolimits_{j}a_{ij}\upsilon _{k;j}=\sum\nolimits_{j}b_{ij}\upsilon
_{k;j}=\lambda _{k}\upsilon _{k;i}\text{.}  \label{VN}
\end{equation}%
On the other hand, if we want that $\vec{\vartheta}_{k}$ to become an
eigenvector of $\mathbf{A}$, with eigenvalue $\Lambda $, then we must impose
the following equations: 
\begin{subequations}
	\label{S}
	\begin{eqnarray}
	\sum_{\ell =1}^{M}a_{i,\ell }\vartheta _{k;\ell }
	&=&\sum\nolimits_{j}b_{i,j}\upsilon _{k;j}+a_{i,M}\mu =\Lambda \upsilon
	_{k;i};\qquad \text{(for any }i\text{)}\label{S1} \\
	\sum_{\ell =1}^{M}a_{M,\ell }\vartheta _{k;\ell }
	&=&\sum\nolimits_{j}a_{M,j}\upsilon _{k,j}+a_{M,M}\mu =\Lambda \mu \text{.}
	\label{S2}
	\end{eqnarray}
	
	Note that in (\ref{S1}) we can use (\ref{VN}) to rewrite the system (\ref{S}%
	) as 
\end{subequations}
\begin{subequations}
	\label{SS}
	\begin{eqnarray}
	a_{iM}\mu &=&\left( \Lambda -\lambda _{k}\right) \upsilon _{k;i};\qquad 
	\text{(for any }i\text{)}  \label{SS1} \\
	\sum\nolimits_{j}a_{Mj}\upsilon _{k;j} &=&\left( \Lambda -a_{MM}\right) \mu 
	\text{.}  \label{SS2}
	\end{eqnarray}
	
	For each one of the $N$ equations in (\ref{SS1}) we conclude that there can
	be only a single solution if each one is multiple of a same equation.
	Therefore, we must impose that the vector $\vec{\nu}$, defined by $\vec{\nu}%
	=\left( a_{1,M},a_{2,M},\ldots ,a_{N,M}\right) ^{\intercal }$, satisfies the
	equality 
\end{subequations}
\begin{equation}
\vec{\nu}=\left( a_{1,M},a_{2,M},\ldots ,a_{N,M}\right) ^{\intercal }=\alpha
\left( \upsilon _{k;1},\upsilon _{k;2},\ldots ,\upsilon _{k;N}\right)
^{\intercal }=\alpha \vec{\upsilon}_{k}\text{,}  \label{VP}
\end{equation}%
such that, assuming $\alpha \neq 0$, we obtain 
\begin{equation}
\mu =\frac{\Lambda -\lambda _{k}}{\alpha }\text{.}  \label{mu}
\end{equation}%
Substituting (\ref{mu}) into the equation (\ref{SS2}) we obtain a second
degree equation, which provides $2$ eigenvalues, namely:%
\begin{equation*}
\left( a_{MM}-\Lambda \right) \left( \lambda _{k}-\Lambda \right) =\alpha
\sum\nolimits_{j}a_{Mj}\upsilon _{k;j}=\sum\nolimits_{j}a_{Mj}a_{jM}\text{,}
\end{equation*}%
where in the last equality, we use (\ref{VP}). Again note that if matrix $%
\mathbf{A}$ is Hermitian, that $\sum_{j=1}^{N}a_{M,j}a_{j,M}=\left\vert
\alpha \right\vert ^{2}\vec{\upsilon}_{k}^{\dag }\cdot \vec{\upsilon}%
_{k}=\left\vert \alpha \right\vert ^{2}$, since we are assuming that the set
of eigenvectors are orthonormal.

The other roots will be determined next, but before let us recall that the
set of eigenvectors $\left\{ \vec{\upsilon}_{k}\right\} $ of $\mathbf{B}$
satisfies the following properties for any $k$ and $j$: 
\begin{subequations}
	\label{P}
	\begin{eqnarray}
	\sum\nolimits_{i}\upsilon _{k;i}^{\ast }\upsilon _{j;i} &=&\delta _{k,j},
	\label{P1} \\
	\sum\nolimits_{i}\upsilon _{k,i}^{\ast }b_{i,j} &=&\lambda _{k}\upsilon
	_{k;j}^{\ast },  \label{P2} \\
	\sum\nolimits_{i}b_{j,i}\upsilon _{k;i} &=&\lambda _{k}\upsilon _{k;j}.
	\label{P3}
	\end{eqnarray}%
	It is also important to remark that since the set of vectors $\left\{ \vec{%
		\upsilon}_{k}\right\} $ form a complete base in an $N$-dimensional space,
	that at least one of the elements of each eigenvector is nonzero. We can
	then label the eigenvectors conveniently so that we can assert that the $i$%
	-th element of the eigenvector $\vec{\upsilon}_{i}$ is non-zero, without
	loss of generality. With these points in mind the demonstration proceeds as
	follows: when calculating $\det \left( \mathbf{A}-\Lambda \mathbf{I}\right) $%
	, we will multiply the line $i$ by $\upsilon _{i;i}^{\ast }$ and so that the
	determinant does not change we will divide it by the same quantity $\upsilon
	_{i;i}^{\ast }$. Next, we must add to the elements of line $i$ the elements
	of each line $j\left( \neq i\right) $ multiplied by $\upsilon _{i;j}^{\ast }$%
	. In this way we obtain, for $i\neq k$, that the determinant $\det \left( \mathbf{A}-\Lambda \mathbf{I}\right)$ becomes
\end{subequations}
\begin{eqnarray}
\frac{\lambda _{i}-\Lambda 
}{\upsilon _{i;i}^{\ast }}\det \left( 
\begin{array}{cccccc}
b_{1,1}-\Lambda & \cdots & b_{1,i} & \cdots & b_{1,N} & \alpha \upsilon
_{k;1} \\ 
\vdots & \ddots &  &  & \vdots & \vdots \\ 
\upsilon _{i;1}^{\ast } &  & \upsilon _{i;i}^{\ast } &  & \upsilon
_{i;N}^{\ast } & 0 \\ 
\vdots &  &  & \ddots & \vdots & \vdots \\ 
b_{N,1} & \cdots & b_{N,i} & \cdots & b_{N,N}-\Lambda & \alpha \upsilon
_{k;N} \\ 
a_{M,1} & \cdots & a_{M,i} & \cdots & a_{M,N} & a_{M,M}-\Lambda%
\end{array}%
\right) \text{.}\notag\\  \label{R0}
\end{eqnarray}%
Since this reasoning holds for all $i\neq k$, it is easy to see that all
eigenvalues $\lambda _{i}$ of $\mathbf{B}$, excluding $i=k$, are also
eigenvalues of $\mathbf{A}$. Therefore, we conclude that%
\begin{eqnarray}
\det \left( \mathbf{A}-\Lambda \mathbf{I}\right)& =&\left[ \left(
a_{M,M}-\Lambda \right) \left( \lambda _{k}-\Lambda \right) -\alpha
\sum\nolimits_{i}a_{M,i}\upsilon _{k;i}\right]\prod\nolimits_{j\left( \neq
	k\right) }\left( \lambda _{j}-\Lambda \right) \text{.}  \label{R1}
\end{eqnarray}

The eigenvectors associated to each one of the eigenvalues common to both
matrices $\mathbf{A}$ and $\mathbf{B}$ can be obtained through the system of
equations:%
\begin{eqnarray*}
	\sum\nolimits_{j}b_{i,j}\upsilon _{k;j}+a_{i,M}\mu  &=&\lambda _{k}\upsilon
	_{k;i}, \\
	\sum\nolimits_{j}a_{M,j}\upsilon _{k;j} &=&\left( \lambda
	_{k}-a_{M,M}\right) \mu \text{,}
\end{eqnarray*}%
from which one realizes that in the case where $\mathbf{A}$ is Hermitian,
such that $\sum\nolimits_{j}a_{M,j}\upsilon _{k;j}\varpropto \vec{\upsilon}%
_{j}^{\dag }\cdot \vec{\upsilon}_{k}=0$ and taking $\mu =0$, we have as the
eigenvector of $\mathbf{A}$, associated with the eigenvalue $\lambda _{j}$, $%
\vec{\vartheta}_{j}=\left( \vec{\upsilon}_{j},0\right) $.

\subsubsection{Step 2: Case in which $M=N+1$ and the elements $\left\{
	a_{i,M}\right\} $ are a linear superposition of all eigenvectors of $\mathbf{%
		B}$}

Now let's look at the situation in which the coefficients of the $M$-th
column of matrix $\mathbf{A}$ - which by definition, compose the vector $%
\vec{\nu}=\left( a_{1,M},a_{2,M},\ldots ,a_{N,M}\right) ^{\top }$ -- are
arbitrary, that is, instead of the relation we have imposed on (\ref{VP}),
let us assume the general case in which the vector $\vec{\nu}$ can always be
written as a linear combination of the form $\vec{\nu}=\sum\nolimits_{j}%
\alpha _{j}\vec{\upsilon}_{j}$, where the coefficients $\alpha _{j}$ are
determined by the inner product $\alpha _{j}=\vec{\upsilon}_{j}^{\dag }\cdot 
\vec{\nu}$.

With the help of equation (\ref{R1}) together with the basic properties of
the determinant, it is easy to verify that the eigenvalues of $\mathbf{A}$
must satisfy the relation%
\begin{eqnarray}
\det \left( \mathbf{A}-\Lambda \mathbf{I}\right) &=&\det \left( 
\begin{array}{cccc}
b_{1,1}-\Lambda & \cdots & b_{1,N} & a_{1,M} \\ 
\vdots & \ddots &  & \vdots \\ 
b_{N,1} &  & b_{N,N}-\Lambda & a_{N,M} \\ 
a_{M,1} & \cdots & a_{M,N} & a_{M,M}-\Lambda%
\end{array}%
\right)=\sum\nolimits_{k}\det \left( 
\begin{array}{cccc}
b_{1,1}-\Lambda & \cdots & b_{1,N} & \alpha _{k}\upsilon _{k;1} \\ 
\vdots & \ddots &  & \vdots \\ 
b_{N,1} &  & b_{N,N}-\Lambda & \alpha _{k}\upsilon _{k;N} \\ 
a_{M,1} & \cdots & a_{M,N} & \frac{a_{M,M}-\Lambda }{N}%
\end{array}%
\right)  \notag \\
&=&\sum\nolimits_{k}\left[ \frac{\left( a_{M,M}-\Lambda \right) \left(
	\lambda _{k}-\Lambda \right) }{N}-\alpha
_{k}\sum\nolimits_{i}a_{M,i}\upsilon _{k;i}\right]\prod\nolimits_{j\left(
	\neq k\right) }\left( \lambda _{j}-\Lambda \right) \text{.}  \label{detG}
\end{eqnarray}

In the situation where $\mathbf{A}$ is Hermitian, so that $%
a_{M,i}=a_{i,M}^{\ast }$, (\ref{detG}) reduces to the form%
\begin{eqnarray}
\det \left( \mathbf{A}-\Lambda \mathbf{I}\right) &=&\sum\nolimits_{k}\left[ 
\frac{\left( a_{M,M}-\Lambda \right) \left( \lambda _{k}-\Lambda \right) }{N}%
-\left\vert \alpha _{k}\right\vert ^{2}\right]\prod\nolimits_{j\left( \neq
	k\right) }\left( \lambda _{j}-\Lambda \right) \text{.}  \label{ca}
\end{eqnarray}

From the expression (\ref{ca}) we find that if one of the eigenvalues $%
\lambda _{k}$ is degenerate, the characteristic polynomial will decrease its
order, since the term $\left( \lambda _{k}-\Lambda \right) $ can be factored
(if the multiplicity of this root is $L\leq N$, then the matrix with
dimension $M=N+1$ will have the same root with the multiplicity $L-1$). The
same occurs when one of the $\alpha _{k}$ is null, since $\left( \lambda
_{k}-\Lambda \right) $ will be present in all terms of this sum. We
therefore have a way of decreasing a degree of the characteristic polynomial
for each imposed restriction, that is, $\alpha _{k}=0$. If we impose that
all $\alpha _{k}$ are null except one, the characteristic polynomial becomes
of order $2$, as in the case of the previous section.

In the next section we will extend the index $\ell =1,\ldots ,L$ to
facilitate the labeling of the columns of matrix $\mathbf{A}$ that are not
part of $\mathbf{B}$.

\subsubsection{Step 3: Case in which $M=N+L$ and the elements of each column $%
	\left\{ a_{i,N+\ell }\right\} $ are a linear superposition of all
	eigenvectors of $\mathbf{B}$}

Now we will verify the influence of the knowledge of eigenvalues for the
case where we have $M=N+L>N+1$. For this situation, let us assume that the
vector $\vec{\nu}_{\ell }$, whose coefficients represent the first $N$
elements of column $N+\ell $ of matrix $\mathbf{A}$, is proportional to a
superposition of all eigenvectors of $\mathbf{B}$, so that we can write $%
\mathbf{A}$ in the form%
\begin{equation*}
\mathbf{A}=\left( 
\begin{array}{cc}
\mathbf{B} & 
\begin{array}{ccc}
\sum_{j}\alpha _{1,j}\upsilon _{j;1} & \cdots  & \sum_{j}\alpha
_{L,j}\upsilon _{j;1} \\ 
\vdots  &  & \vdots  \\ 
\sum_{j}\alpha _{1,j}\upsilon _{j;N} & \cdots  & \sum_{j}\alpha
_{L,j}\upsilon _{j;N}%
\end{array}
\\ 
\begin{array}{ccc}
a_{N+1,1} & \cdots  & a_{N+1,N} \\ 
\vdots  &  & \vdots  \\ 
a_{N+L,1} & \cdots  & a_{N+L,N}%
\end{array}
& 
\begin{array}{ccc}
a_{N+1,N+1} & \cdots  & a_{N+1,N+L} \\ 
\vdots  & \ddots  & \vdots  \\ 
a_{N+L,N+1} & \cdots  & a_{N+L,N+L}%
\end{array}%
\end{array}%
\right) \text{.}
\end{equation*}

Using (\ref{P1}) and (\ref{P2}) and following the same steps we perform to
obtain (\ref{R0}), we can show that $\det \left( \mathbf{A}-\Lambda \mathbf{I%
}\right) $ can be written as%
\begin{equation}
\det \left( 
\begin{array}{cccccc}
b_{1,1}-\Lambda  & \cdots  & b_{1,N} & \sum_{j}\alpha _{1,j}\upsilon _{j;1}
& \cdots  & \sum_{j}\alpha _{L,j}\upsilon _{j;1} \\ 
\vdots  &  & \vdots  & \vdots  &  & \vdots  \\ 
\left( \lambda _{i}-\Lambda \right) \frac{\upsilon _{i;1}^{\ast }}{\upsilon
	_{i;i}^{\ast }} & \cdots  & \left( \lambda _{i}-\Lambda \right) \frac{%
	\upsilon _{i;N}^{\ast }}{\upsilon _{i;i}^{\ast }} & \frac{\alpha _{1,i}}{%
	\upsilon _{i;i}^{\ast }} & \cdots  & \frac{\alpha _{L,i}}{\upsilon
	_{i;i}^{\ast }} \\ 
\vdots  &  &  & \vdots  &  & \vdots  \\ 
b_{N,1} & \cdots  & b_{N,N}-\Lambda  & \sum_{j}\alpha _{1,j}\upsilon _{j;N}
& \cdots  & \sum_{j}\alpha _{L,j}\upsilon _{j;N} \\ 
a_{N+1,1} & \cdots  & a_{N+1,N} & a_{N+1,N+1}-\Lambda  & \cdots  & 
a_{N+1,N+L} \\ 
\vdots  &  & \vdots  & \vdots  & \ddots  & \vdots  \\ 
a_{N+L,1} & \cdots  & a_{N+L,N} & a_{N+L,N+1} & \cdots  & 
a_{N+L,N+L}-\Lambda 
\end{array}%
\right) \text{.}  \label{PP}
\end{equation}

From (\ref{PP}), we note that if none of the $\vec{\nu}_{\ell }$ is
proportional to the eigenvector $\vec{\upsilon}_{i}$, such that $\alpha
_{\ell ,i}=0$, the eigenvalue $\lambda _{i}$ of $\mathbf{B}$ will also be a
eigenvalue of matrix $\mathbf{A}$, thus reducing the order of the
characteristic polynomial. However, note that the minimum order of the
characteristic polynomial in this case will be $M-\left( N-1\right)
=N+L-N+1=L+1$, which corresponds to the case where each one of $\vec{\nu}%
_{\ell }$ is proportional to the same eigenvector of $\mathbf{B}$, i.e. $%
\vec{\nu}_{\ell }=\beta _{\ell }\vec{\upsilon}_{i}$.

An expression analogous to (\ref{detG}) here is not possible, except for $%
\vec{\nu}_{\ell }=\beta _{\ell }\vec{\upsilon}_{i}$, because we obtain a
system of equations that can not be solved without imposing constraints on
the elements of matrix $\mathbf{B}$. However, we have seen in the previous
section that if a $N\times N$ square matrix has a root with multiplicity $%
L\leq N$, then the matrix with dimension $M=N+1$, regardless of the elements
of column $N+1$, will have the same root with a multiplicity of at least $L-1
$. Following with this reasoning, by induction, we conclude that a matrix
with dimension $M=N+2$, with any elements of the columns $N+1$ and $N+2$, it
will have the same root with at least one multiplicity $L-2$ and so on, such
that, the matrix with dimension $M=N+L-1$ will have the same root with at
least one multiplicity $1$.

Finally, it is important to mention that the properties we develop here for
the columns can also be applied to the rows, since the eigenvalues of a
matrix do not change if we take its transpose. Thus, if there are no
satisfied constraints for the columns, one must check if they become valid
for the rows when the matrix under analysis is not Hermitian.

\section{Applications}

In this section, we try to apply the formalism developed above in two
situations of interest: $i)$ to expand a subset of eigenvectors in a space $%
%TCIMACRO{\U{2102} }%
%BeginExpansion
\mathbb{C}
%EndExpansion
^{N}$ to a space of greater dimension preserving the modulus, direction and
sense of these eigenvectors, and $ii)$ to construct an square matrices with
arbitrarily dimensioned, in which their eigenvalues can be determined
analytically.

\subsection{Bringing a set of eigenvectors from the space $%
	%TCIMACRO{\U{2102} }%
	%BeginExpansion
	\mathbb{C}
	%EndExpansion
	^{N}$ to the space $%
	%TCIMACRO{\U{2102} }%
	%BeginExpansion
	\mathbb{C}
	%EndExpansion
	^{N+L}$, preserving its modulus, direction and sense}

Consider a square matrix $\mathbf{B}$, defined by%
\begin{equation}
\mathbf{B}_{N\times N}=\left( 
\begin{array}{cccc}
b_{1,1} & b_{1,2} & \cdots  & b_{1,N} \\ 
b_{2,1} & b_{2,2} & \cdots  & b_{2,N} \\ 
\vdots  & \vdots  & \ddots  & \vdots  \\ 
b_{N,1} & b_{N,2} & \cdots  & b_{N,N}%
\end{array}%
\right) ,  \label{Q}
\end{equation}%
in which its eigenvalues $\left\{ \lambda _{j}\right\} $ and respective
eigenvectors $\left\{ \vec{\upsilon}_{j}\right\} $ are known. According to
the theory above, we can extend the dimension of this matrix to $N+1$,
keeping up to $N-1$ eigenvectors of $\mathbf{B}$, extending them to the
space $%
%TCIMACRO{\U{2102} }%
%BeginExpansion
\mathbb{C}
%EndExpansion
^{N+1}$ adding a null component. As an example, let's say we want to ``keep"
the set of eigenvectors $\left\{ \vec{\upsilon}_{j\neq k}\right\} $ of $%
\mathbf{B}_{N\times N}$. To perform this procedure, we must choose a
partially Hermitian matrix of the form 
\begin{equation*}
\mathbf{B}_{\left( N+1\right) \times \left( N+1\right) }=\left( 
\begin{array}{cc}
\mathbf{B} & \alpha \vec{\upsilon}_{k} \\ 
\beta \vec{\upsilon}_{k}^{\dagger } & a_{N+1,N+1}%
\end{array}%
\right) \text{.}
\end{equation*}

This procedure can be applied countless times in order to obtain a larger
dimension matrix that ``preserves" the eigenvectors of $\mathbf{B}_{N\times N}
$. One of the possibilities becomes:%
\begin{equation}
\mathbf{B}_{\left( N+L\right) \times \left( N+L\right) }=\left( 
\begin{array}{ccccc}
\mathbf{B} & \alpha _{1}\vec{\upsilon}_{k} & \alpha _{2}\vec{\upsilon}_{k} & 
\cdots  & \alpha _{L}\vec{\upsilon}_{k} \\ 
\beta _{1}\vec{\upsilon}_{k}^{\dagger } & a_{N+1,N+1} & \alpha _{2}\frac{%
	\Lambda _{1}-\lambda _{k}}{\alpha _{1}} & \cdots  & \alpha _{L}\frac{\Lambda
	_{1}-\lambda _{k}}{\alpha _{1}} \\ 
\beta _{2}\vec{\upsilon}_{k}^{\dagger } & \beta _{2}\left( \frac{\Lambda
	_{1}-\lambda _{k}}{\alpha }\right) ^{\ast } & a_{N+2,N+2} & \cdots  & \alpha
_{L}\frac{\Lambda _{2}-\lambda _{k}}{\alpha _{2}} \\ 
\vdots  & \vdots  & \vdots  & \ddots  & \vdots  \\ 
\beta _{L}\vec{\upsilon}_{k}^{\dagger } & \beta _{L}\left( \frac{\Lambda
	_{1}-\lambda _{k}}{\alpha }\right) ^{\ast } & \beta _{L}\left( \frac{\Lambda
	_{2}-\lambda _{k}}{\alpha _{2}}\right) ^{\ast } & \cdots  & a_{N+L,N+L}%
\end{array}%
\right) \text{,}  \label{O}
\end{equation}%
considering $\Lambda _{0}=\lambda _{k}$ and $\Lambda _{\ell }$ ($\ell
=1,2,\ldots ,L$) being one of the roots of the 2nd degree polynomial%
\begin{equation*}
\left( \Lambda _{\ell }-a_{N+\ell ,N+\ell }\right) \left( \Lambda _{\ell
}-\Lambda _{\ell -1}\right) -\beta _{\ell }^{\ast }\alpha _{\ell }=0.
\end{equation*}%
The other possibilities arise from the fact that in $\mathbf{B}_{\left(
	N+2\right) \times \left( N+2\right) }$ we can form a column by the linear
combination of the two ``excluded" eigenvectors of the space we wish to
preserve and not a single one as we did in (\ref{O}). Similarly in $\mathbf{B%
}_{\left( N+3\right) \times \left( N+3\right) }$ we can form a column by the
linear combination of 3 of the excluded eigenvectors and so on.

It is important to note that the elements $a_{N+\ell ,N+\ell }$ of $\mathbf{B%
}_{\left( N+L\right) \times \left( N+L\right) }$, for $\ell $ ranging from $1
$ to $L$, are only important for determining the eigenvalues other than $%
\left\{ \lambda _{j\neq k}\right\} $, which does not interfere in the set of
eigenvectors $\left\{ \vec{\upsilon}_{j\neq k}\right\} $, which are extended
into the vector space defined by the matrix $\mathbf{B}_{\left( N+L\right)
	\times \left( N+L\right) }$. Note that this property in quantum mechanics
can be very interesting because it shows us what are the necessary bonds
that must be imposed on the quantum system in order for a set of self-states
to be ``preserved" or taken to a larger Hilbert space, such that we can
guarantee some particular solutions in this larger Hilbert space.

In this context, note that this paper provides a protocol to make
approximations in the numerical calculation of eigenvalues where the
coefficients representing the constraints, although not be satisfied, can be
considered to be approximately valid when they are much smaller than the
coefficients of the other eigenvectors.

\subsection{Identifying arbitrary $N\times N$ square matrices in which its
	eigenvalues can be determined analytically}

In the literature we find some matrices with certain symmetries that make it
possible to obtain the eigenvalues and eigenvectors for cases in which the
dimension of the square matrix is arbitrary, as we mentioned in the
introduction. With the study we have developed here, we can significantly
expand the number of situations in which we can analytically resolve such
matrices, as we see below.

The process is very similar to that of the previous section, but here there
is no need to maintain an eigenvectors subspace when we increase the size of
the matrix. This process can be seen in an inverse way, that is, given the
final matrix, we impose constraints on some elements in order to make
possible the calculation of eigenvalues and eigenvectors algebraically.

Consider the same square matrix (\ref{Q}) above. Let us assume that this $%
N\times N$ square matrix represents any matrix that we know how to calculate
the eigenvalues analytically, and may even be one of the $4$ analytical
matrices mentioned in the introduction. Our desire is to obtain a new square
matrix of larger size, in which the eigenvalues can be found in an
analytical way, and therefore, that its characteristic polynomial can be
decomposed into products of polynomials of, at most, order $4$. For this
reason, in constructing the $N+1$ dimension matrix we should at least
restrict the elements of the $N+1$ column to a linear combination of $3$
eigenvectors of the matrix $\mathbf{B}_{N\times N}$ so that the new matrix $%
\mathbf{B}_{\left( N+1\right) \times \left( N+1\right) }$ has $N-3$
eigenvalues common to $\mathbf{B}_{N\times N}$, while the remaining
eigenvalues ($4$ remain) can be obtained from the root of a fourth-order
polynomial. This process can be repeated successively until a desired high
size matrix is obtained. After the first development the eigenvalues will be
the following: we will have $N-4$ eigenvalues identical to the matrix $%
\mathbf{B}_{N\times N}$ that we will order as $\lambda _{1},\lambda
_{2},\ldots ,\lambda _{N-3}$. The other eigenvalues will be obtained from
the root of the 4th order polynomial%
\begin{eqnarray*}
	P\left( \Lambda \right)  &=&\left( a_{N+1,N+1}-\Lambda \right) \left(
	\lambda _{N}-\Lambda \right) \left( \lambda _{N-1}-\Lambda \right) \left(
	\lambda _{N-2}-\Lambda \right)  \\
	&&-\left\vert \alpha _{N-2}\right\vert ^{2}\left( \lambda _{N}-\Lambda
	\right) \left( \lambda _{N-1}-\Lambda \right) -\left\vert \alpha
	_{N-1}\right\vert ^{2}\left( \lambda _{N}-\Lambda \right) \left( \lambda
	_{N-2}-\Lambda \right)-\left\vert \alpha _{N}\right\vert ^{2}\left( \lambda _{N-1}-\Lambda
	\right) \left( \lambda _{N-2}-\Lambda \right)  \\
	&=&\det \left( 
	\begin{array}{cccc}
		\lambda _{N}-\Lambda  & 0 & 0 & \left\vert \alpha _{N}\right\vert  \\ 
		0 & \lambda _{N-1}-\Lambda  & 0 & \left\vert \alpha _{N-1}\right\vert  \\ 
		0 & 0 & \lambda _{N-2}-\Lambda  & \left\vert \alpha _{N-2}\right\vert  \\ 
		\left\vert \alpha _{N}\right\vert  & \left\vert \alpha _{N-1}\right\vert  & 
		\left\vert \alpha _{N-2}\right\vert  & a_{N+1,N+1}-\Lambda 
	\end{array}%
	\right) .
\end{eqnarray*}%
By labeling these roots as $\lambda _{N-2}^{\left( 1\right) },\lambda
_{N-1}^{\left( 1\right) },\lambda _{N}^{\left( 1\right) }$ and $\lambda
_{N+1}^{\left( 1\right) }$ and the others as $\lambda _{k}^{\left( 1\right)
}=\lambda _{k}$ for $k=1,2,\ldots ,N-3$ we return to the initial problem of
the first development, replacing $N$ with $N^{\left( 1\right) }=N+1$. In the
second development we can choose the elements of the $N+2$ column as a
linear combination of any of the $3$ eigenvectors of the matrix $\mathbf{B}%
_{\left( N+1\right) \times \left( N+1\right) }$ and so on. After $L$
processes, we obtain the following matrix%
\begin{equation*}
\mathbf{B}_{\left( N+L\right) \times \left( N+L\right) }=\left( 
\begin{array}{cc}
\begin{array}{cc}
\begin{array}{cc}
\mathbf{B} & \sum\limits_{j=N-2}^{N}\alpha _{j}\vec{\upsilon}_{j}^{\left(
	N\right) } \\ 
c.t. & a_{N+1,N+1}%
\end{array}
& \cdots  \\ 
\vdots  & \ddots 
\end{array}
& \sum\limits_{j=N+L-2}^{N+L}\alpha _{j}\vec{\upsilon}_{j}^{\left(
	N+L-1\right) } \\ 
c.t. & a_{N+L,N+L}%
\end{array}%
\right) ,
\end{equation*}%
where in each line with the letters ``$c.t.$" means the conjugate transpose
of the respective column, such that the matrix $\mathbf{B}_{\left(
	N+L\right) \times \left( N+L\right) }$ is Hermitian.

Note that this technique can be used for alternative cryptographic protocols
where matrix $\mathbf{B}$ can be shared by a public key and the sequence of
eigenvectors on a private channel or vice versa.

\section{Summary and conclusions}

In this work we show that the knowledge of the eigenvalues and eigenvectors
of a matrix $\mathbf{B}_{N\times N}$, which is a submatrix of $\mathbf{A}%
_{M\times M}$ obtained by deleting the respective rows and columns, may help
us to reduce the degree of the characteristic polynomial associated with
matrix $\mathbf{A}$. This reduction can come in two ways: $i)$ by imposing
links of a some elements of $\mathbf{A}$ and $ii)$ when any of the
eigenvalues of $\mathbf{B}$ have degeneracy greater than $M-N$. The latter
case exempts any imposition of restriction and is valid for any square
matrix, not necessarily Hermitian. This results are highlighted in $2$
simple applications. In the first application, we show how to extend an
eigenvector in $%
%TCIMACRO{\U{2102} }%
%BeginExpansion
\mathbb{C}
%EndExpansion
^{N}$ to a space $%
%TCIMACRO{\U{2102} }%
%BeginExpansion
\mathbb{C}
%EndExpansion
^{N+L}$, preserving its norm, direction and sense. This application can be
used in quantum system studies in which it is desired that an interacting
system (associated with a Hilbert space of dimension $N$) preserves its
self-states when we include new interactions (in which we can associate a
Hilbert space of dimension $N+L$). In the second application, we show how to
construct matrices of arbitrary dimensions, nontrivial, that allow
analytical calculations for eigenvalues and that can be used in
cryptographic protocols.

In order to finish this work, whose objective was only to illustrate our
technique, apparently absent in the literature, we mention some of the
possible horizons that emerge from this work and that need an additional
studies:

\begin{enumerate}
	\item Formalize a protocol or algorithm to determine the eigenvalues taking
	into account the technique developed in this work in conjunction with those
	already in the literature and to verify if there is a computational time
	gain;
	
	\item Formalize a protocol or algorithm in cryptography based on the
	properties of the matrices developed in this work;
	
	\item Using a random selection of eigenvectors in the application example 2,
	investigate whether these analytic matrices can also be considered random or
	have a small degree of order;
	
	\item Verify if the dynamics between states of a subspace of dimension $N$,
	whose eigenvectors are extended to another of dimension $N+L$, remain intact
	for some physical situations of interest.
	
	\item What modifications should we make in this work to consider a base of
	eigenvectors of matrix $\mathbf{B}$ not orthonormal? Is this possible?
\end{enumerate}

\begin{acknowledgments}
L.C.C. would like to thank FAPESP for the financial support during its IC
research and M. A. de Ponte would like to thank G.I.Y. by encouragement and
motivation and my former university professors A.K.M. Libardi and I.C.
Rossini and my elementary school teacher Midouri.
\end{acknowledgments}

\end{document}